\documentclass[12pt]{amsart}
\usepackage[a4paper,margin=23mm]{geometry}

\usepackage{amsmath, amssymb}
\usepackage{amsthm}
\usepackage{graphicx}
\usepackage{comment}
\usepackage[colorlinks=true, allcolors=blue]{hyperref}

\newtheorem{theorem}{Theorem}[section]

\newtheorem{remark}[theorem]{Remark} 
\newtheorem{lemma}[theorem]{Lemma}
\usepackage[english]{babel}

\usepackage[backend=bibtex,
style=numeric,
isbn=false,
doi=false
]{biblatex}
\title{Binomial coefficients coprime to 6}
\author{Pascal Jelinek}
\email{pascal.jelinek@unileoben.ac.at}
\address{
	Department Mathematics and Information Technology,
	Leoben University of Technology,
	Franz-Josef-Strasse 18, 8700 Leoben, Austria.\
}

\begin{document}
\maketitle
\renewcommand{\thefootnote}{\fnsymbol{footnote}}
\footnotetext{\emph{2020 Mathematics Subject Classification.} Primary: 11B65, 05A10; Secondary: 11A07,11B50}
%\footnotetext{\emph{Key words and phrases.}  binomial coefficient, }
\footnotetext{ The author was supported by the Austrian Science Fund (FWF), Project P36137-N.}
\renewcommand{\thefootnote}{\arabic{footnote}}
\begin{abstract}
	It is a well known result by Singmaster that any integer $d$ divides almost all binomial coefficients. The study of the structure and the size of the exceptional set has been of interest for many decades. In the case where $d$ is a prime power, these sets are well understood and their size is known. However, if $d$ has at least two distinct prime factors, no non-trivial bounds are known. In this paper, we will provide the first non trivial bound on the number of binomial coefficients coprime to 6.
\end{abstract}

\section{Introduction}
The divisibility of binomial coefficients has been studied extensively for a long time. In 1974, Singmaster \cite{Singmaster1974} showed that every natural number $d$ divides almost all binomial coefficients. More precisely, his result states that as $N\rightarrow\infty$, we have for all $d$ that
\[
\frac{2}{N(N-1)}\#\left\{0\leq m\leq n< N: d\,\middle\vert\,\binom{n}{m}\right\}\rightarrow 1\,.
\]
Therefore, we are interested in studying the size and structure of the exceptional set.

For $d=p$, Lucas \cite{Lucas1878} proved the following theorem:
\begin{theorem}[Lucas]\label{Lucasthm}
	Let $p$ be a prime number. Let $m\leq n$ be two non-negative integers. Consider the base $p$ expansions of $n=\sum\limits_{i=0}^{k}n_ip^i$ and of $m=\sum\limits_{i=0}^{k}m_ip^i$. Then
	\begin{align}
		\binom{n}{m}\equiv \prod_{i=0}^{k} \binom{n_i}{m_i} \bmod p\,.
	\end{align}
\end{theorem}
Here we use the convention that $\binom{n}{m}=0$, if $m>n$. Hence, $p \nmid \binom{n}{m}$ if and only if $n_i\geq m_i$ for all $i$. 
In the case $p=2$, Stolarsky \cite{Stolarsky1977} and Harborth \cite{Harborth1976} even proved that there exists a continuous, nowhere differentiable, 1-periodic function c(n) such that for all $N>1$
\begin{align}
	\#\left\{0\leq m\leq n < N: 2\not\,\middle\vert\, \binom{n}{m}\right\} = c\left(\frac{N}{2^{\lceil \log_2(N)\rceil}}\right)\cdot N^{\frac{\log(3)}{\log(2)}}.
\end{align}
This was later generalized by Stein \cite{Stein1989}.
\begin{theorem} [Stein]\label{thm_Stein}
	Let $p$ be a prime number. Then there exists a continuous, nowhere differentiable, 1-periodic function c(n) such that for all $N>1$
	\begin{align}
		\#\left\{0\leq m\leq n < N: p\not \,\middle\vert\, \binom{n}{m}\right\} = c\left(\frac{N}{p^{\lceil \log_p(N)\rceil}}\right)\cdot N^{\frac{\log(p(p+1)/2)}{\log(p)}}.
	\end{align}
\end{theorem}
Further, it is known by Wilson \cite{Wilson1998} that $1/2<c(n)\leq 1$ and $c(1)=1$.

Theorem \ref{Lucasthm} was independently generalized by Davis and Webb \cite{DavisWebb1990}, and by Granville \cite{Granville1997} to the case where $d$ is a prime power, say $d=p^a$. Barat and Grabner \cite{BaratGrabner2001} used the formulation of Granville to show that all binomial coefficients with the same $p$-valuation are equidistributed over the possible residue classes mod $p^a$. Spiegelhofer and Wallner \cite{SpiegelhoferWallner2023} showed that the $p$-valuation of binomial coefficients is normally distributed.

However, to the author's knowledge, no upper bounds on the size of the exceptional set are known in the case where d is not a prime power. We want to study the case $d=6$ in this paper.

Using the result by Singmaster combined with the Chinese Remainder Theorem, one can immediately prove the size of the sets for the residue classes $0,2,3,4$. We have that, there exit constants $c_2, c_3$ such that
\begin{align}
	\frac{2}{N(N-1)}&\#\left\{0\leq m\leq n< N: 6\,\middle\vert\,\binom{n}{m}\equiv 0 \right\}\rightarrow 1\\
	c_3 N^{\log6/\log3}\leq &\#\left\{0\leq m\leq n < N:\gcd\left( \binom{n}{m},6\right) = 2 \right\} \leq N^{\log6/\log3}\\
	c_2  N^{\log3/\log2} \leq &\#\left\{0\leq m\leq n < N: \gcd\left(\binom{n}{m},6\right) = 3\right\} \leq N^{\log3/\log2}\; .
\end{align}
For the residue classes $1$ and $5$, we heuristically expect by the independence of primes that
\[ 
\#\left\{0\leq m\leq n < N: \gcd\left(\binom{n}{m},6\right) = 1 \right\} \ll N^{\log3/\log2+\log6/\log3-2}.
\]
On the other hand, one can bound these sets trivially by the number of odd binomial coefficients and therefore get the following inequality: 
\[
\#\left\{0\leq m\leq n < N: \gcd\left(\binom{n}{m},6\right) = 1 \right\} \leq N^{\log3/\log2}.
\]

In this paper we will provide log-savings for the trivial bound. Namely, we get the following theorem:
\begin{theorem}\label{thm_new}
	Let $N>0$ be an integer. Then
	\begin{equation}
		\#\left\{0\leq m\leq n < N: \gcd\left(\binom{n}{m},6\right) = 1 \right\}
		\ll \frac{N^{\log3/\log2}}{\log(N)^{2(1-\log_3(2))-\varepsilon}}.
	\end{equation}
\end{theorem}
\begin{comment}
	\section{Acknowledgements}
	The author wants to thank Lukas Spiegelhofer for introducing him to this problem. The author was supported by the FWF project P36137.
\end{comment}

\section{Proof of the theorem}
The idea of the proof is to study the number of odd binomial coefficients $\binom{n}{m}$, where we restrict $n$ and $m$ to certain residue classes modulo $3^k$. Since these residue classes correspond to the last $k$ digits of the base $3$ expansion of $n$ and $m$, we will use Stein's theorem (Theorem \ref{thm_Stein}) for $p=3$ to rule out some possible combinations of residue classes of $n$ and $m$, hence generating a saving over the trivial bound, which is just the number of odd binomial coefficients.

Before we provide the details of the proof, we fix some notation. Our binomial coefficients will always be denoted by $\binom{n}{m}$. Further, let $(n_i)_{i\geq0}$, $(m_i)_{i\geq0}$, be the digits of $n$ and $m$, respectively, in the base 2 expansions, i.e.
\[
n=\sum\limits_{i\geq0}n_i2^i
\]
and 
\[
m=\sum\limits_{i\geq0}m_i2^i.
\]
Additionally, let 
\[
\zeta=\exp(2\pi i/3^k)
\] be a $3^k$-th root of unity.

Let $L$ be such that
\[
1<  \frac{2^{L}}{N} \leq 2\,.
\]
Now we begin the proof of Theorem \ref{thm_new}. Let $n<2^L$ be an integer. First we analyse, how many $\binom{n}{m}$ are coprime to $2$. By Lucas's theorem (Theorem \ref{Lucasthm}), we see that whenever $n_i=0$, then $m_i=0$, and whenever $n_i=1$, then $m_i\in \{0,1\}$. Therefore, we see that
\[
\#\left\{0\leq m\leq n: \gcd\left(\binom{n}{m},2\right) = 1 \right\} = 2^{s_2(n)}.
\]
Since 
\[
\sum_{m_i=0}^{n_i}\binom{n_i}{m_i}=2^{n_i},
\]
we can rewrite the above equality as
\[
2^{s_2(n)}=\prod_{i=0}^{{L-1}}2^{n_i}=\prod_{i=0}^{L-1}\sum_{m_i=0}^{n_i}\binom{n_i}{m_i}.
\]
If we sum over $n<2^L$, we see that the number of binomial coefficients that are coprime to $2$ can be expressed as
\begin{equation}\label{eq_startpoint}
	\sum_{n_0=0}^{1}\cdots\sum_{n_{L-1}=0}^{1}\prod_{i=0}^{L-1}\sum_{m_i=0}^{n_i}\binom{n_i}{m_i}.
\end{equation}

So far we have only used the digit expansion in base $2$. We note that we can detect the last $k$ digits of an integer $n$ in base $3$, by evaluating $n \bmod 3^k$. By Stein's theorem (Theorem \ref{thm_Stein}) we see that only $6^k$ out of the possible $3^{2k}$ combinations of the last $k$ digits of $n,m$ are admissible such that $\binom{n}{m}$ can still be coprime to $3$. We now aim to choose $k$ as large as possible with respect to $L$ in order to maximize the savings. 

By the observation above, we now want to study the number of binomial coefficients $\binom{n}{m}$ that are coprime to $2$ and satisfy $n\equiv r_1 \bmod 3^k$ and $m\equiv r_2 \bmod 3^k$. By equation \eqref{eq_startpoint}, this number is
\[
\sum_{n_0=0}^{1}\cdots\sum_{\substack{n_{L-1}=0\\ n\equiv r_1\bmod3^k}}^{1}
\sum_{m_0=0}^{n_i}\cdots\sum_{\substack{m_{L-1}=0\\ m\equiv r_2\bmod3^k}}^{n_{L-1}}\binom{n_0}{m_0}\cdots\binom{n_{L-1}}{m_{L-1}}.
\]
We now use the orthogonality of characters to deduce that the above expression is equal to
\begin{align*}
	\frac 1{3^{2k}}\sum_{a,b=0}^{3^k-1}\zeta^{-ar_1}\zeta^{-br_2}
	&\sum_{n_0}^{1}\cdots\sum_{n_{L-1}}^{1} \zeta^{an_02^0}\cdots\zeta^{an_{L-1}2^{L-1}}\\
	&\sum_{m_0}^{n_0} \binom{n_0}{m_0} \zeta^{bm_02^0} \cdots\sum_{m_{L-1}=0}^{n_{L-1}} \binom{n_{L-1}}{m_{L-1}} \zeta^{bm_{L-1}2^{L-1}}.
\end{align*}
By the binomial theorem, we obtain
\begin{align*}
	\frac 1{3^{2k}}\sum_{a,b=0}^{3^k-1}\zeta^{-ar_1}\zeta^{-br_2}
	&\sum_{n_0}^{1} \zeta^{an_02^0} \left(1+\zeta^{b2^0}\right)^{n_0}\cdots\sum_{n_{L-1}}^{1}\zeta^{an_{L-1}2^{L-1}}\left(1+\zeta^{b2^{L-1}}\right)^{n_{L-1}}.
\end{align*}
Applying the binomial theorem again we get
\[
\frac 1{3^{2k}}\sum_{a=0}^{3^k-1}\sum_{b=0}^{3^k-1}\zeta^{-ar_1}\zeta^{-br_2} \prod_{i=0}^{L-1}\left(1+\zeta^{a2^i}+\zeta^{(a+b)2^i}\right).
\]
If $a=b=0$, we get the main term 
\[
\frac {3^L}{3^{2k}}.
\]
Therefore, it remains to show that
\[
\frac 1{3^{2k}}\sum_{a=0}^{3^k-1}\sum_{\substack{b=0\\(a,b)\neq (0,0)}}^{3^k-1}\zeta^{-ar_1}\zeta^{-br_2} \prod_{i=0}^{L-1}\left(1+\zeta^{a2^i}+\zeta^{(a+b)2^i}\right) = o\left(	\frac {3^L}{3^{2k}}\right).
\]
After taking absolute values and using the triangle inequality on the left-hand side, we aim to show that
\[
\frac 1{3^{2k}}\sum_{a=0}^{3^k-1}\sum_{\substack{b=0\\(a,b)\neq (0,0)}}^{3^k-1} \prod_{i=0}^{L-1}\left|1+\zeta^{a2^i}+\zeta^{(a+b)2^i}\right| = o\left(	\frac {3^L}{3^{2k}}\right).
\]
In order to show that the product is small, we need to analyse the set of differences
\[
\{2^x-2^y\bmod3^k:0\leq x,y\leq L\}.
\]
To achieve this, we will use a special case of a lemma in Korobov \cite[Lemma 2]{Korobov1972}:
\begin{lemma}\label{lem_Kor}
	Let $\gcd(a,3^k)\leq 3^{k-2}$ and $L\leq 2\cdot3^{k-1}$. Then
	\[
	\left|\sum_{i=0}^{L-1} \zeta^{a2^i}\right|\ll \sqrt{3^k}k.
	\]
\end{lemma}
We get a non-trivial estimate of the sum when $L$ is at least $3^{k/2+\varepsilon}$. We will indeed choose $L=3^{k/2+\varepsilon}$.
\begin{remark}
	In the easier case of studying exponential sums of this kind with prime modules, the best currently known relation which is due to Konyagin and Shparlinski \cite{KonyaginShparlinski2012} would only slightly improve the relationship between $L$ and $k$ to
	\[
	L = 3^{12k/25 + \varepsilon}.
	\]
\end{remark}
Using Lemma \ref{lem_Kor}, we get by the Erd\H{o}s-Turan inequality that for any integers $y,H$
\[
\#\{0\leq i \leq L-1: a2^i\in\{y+1,y+2, \dots, y+H\}\bmod3^k\}=\frac{LH}{3^k}+O(E)\;,
\]
where the error term is bounded by
\[
E\leq \frac Ln + \sum_{x=1}^{n}\frac 1x \left|\sum_{i=0}^{L-1} \zeta^{ax2^i}\right|
\]
for any $n$. Taking $n=3^k$, yields
\[
E \leq \frac L{3^k} + \sum_{x=1}^{3^k}\frac 1x \sqrt{3^k}k \leq O\left(\sqrt{3^k}k\log\left(3^k\right)\right) = O\left(\sqrt{3^k}k^2\right).
\]
Next, we take $y=H=3^{k-1}$ and obtain
\[
\#\{0\leq i \leq L-1: a2^i\in\{3^{k-1}+1,3^{k-1}+2, \dots, 2\cdot3^{k-1}\}\bmod3^k\}=\frac{L}{3}+O\left(\sqrt{3^k}k^2\right)
>\frac{L}{4}.
\]
Therefore, for at least $\frac L4$ values of $i$, either $a2^i \bmod 3^k$ or $(a+b)2^i \bmod 3^k$ lies in the interval \[3^{k-1},\dots, 2\cdot3^{k-1}\] and we can bound
\[
\left|1+\zeta^{a2^i}+\zeta^{(a+b)2^i}\right|
\]
from above by
\[
\sqrt{7}.
\]
We bound the remaining factors trivially by $3$ and get that the error term is bounded above by
\[
\frac 1{3^{2k}}\sum_{a=0}^{3^k-1}\sum_{\substack{b=0\\(a,b)\neq (0,0)}}^{3^k-1} \prod_{i=0}^{L-1}\left|1+\zeta^{a2^i}+\zeta^{(a+b)2^i}\right| \leq \sqrt{7}^{\frac L4}3^{\frac{3L}4}=o\left(3^{0.99L}\right).
\]

This concludes our analysis of the number of odd binomial coefficients in an arbitrary residue class modulo $3^k$. We get that the number of odd binomial coefficients in one pair of residue classes is
\[
\frac{3^L}{3^{2k}} +o\left(3^{0.99L}\right).
\]
By Stein's theorem (Theorem \ref{thm_Stein}), there are $6^k$ admissible combinations of residue classes $r_1,r_2$, such that it is possible that $\binom{n}{m}$ is not divisible by $3$. We can bound the number of binomial coefficients $\binom{n}{m}$ that are coprime to $6$ from above by
\[
\ll \frac{6^k3^{L}}{3^{2k}} = \frac{3^L}{\left(\frac{3}{2}\right)^k}.
\]
Therefore, recalling that $2^{L-1}<N\leq2^{L}$ is bounded, we get
\[
\ll \frac{N^{\log3/\log2}}{\left(\frac{3}{2}\right)^k}.
\]
Hence, by using the relation between $N$, $L$ and $k$, we can express $\left(\frac{3}{2}\right)^k$, in terms of $N$ and get that at most
\[
\ll \frac{N^{\log3/\log2}}{\log(N)^{2(1-\log_3(2))-\varepsilon}}
\]
many binomial coefficients are coprime to $6$, hence proving the Theorem \ref{thm_new}.

\begin{remark}
	The author strongly believes that a similar analysis can be done in general for $d=pq$. However, since the saving via this method is still far away from the conjectured exponent these cases are left out.
\end{remark}

\printbibliography
\end{document}